\documentclass[12pt]{article}
 \usepackage{graphicx}
 \usepackage{amssymb}
 \usepackage{epstopdf}
 \newtheorem{thm}{{\bf Theorem}}[section]
  \newtheorem{lem}[thm]{{\bf Lemma}}
  \newtheorem{cor}[thm]{{\bf Corollary}}

\title{The forcing partial order on a family of braids  forced by  pseudo--Anosov $3$--braids}
\author{ Eiko Kin}
\date{}

\begin{document}
 \maketitle
\begin{abstract}
Li--York theorem tells us that  a period $3$ orbit for a continuous map of the interval into itself implies the existence of a periodic orbit of every period. 
This paper concerns an analogue of the theorem for  homeomorphisms of the $2$--dimensional disk. 
In this case a periodic orbit  is specified by a braid type and on the set of all braid types Boyland's dynamical  partial order can be defined. 
We describe the partial order on a family of braids and show that a period $3$ orbit of pseudo--Anosov braid type implies the Smale--horseshoe map 
which is a factor possessing complicated chaotic dynamics. 
\end{abstract}
\footnotetext{

Date: October 18, 2007}
\section{Introduction}

Let $f: D \rightarrow D$ be an  orientation preserving homeomorphism of the $2$-dimensional closed disk $D$. 
One of the main question on dynamical systems is to investigate the variety of periodic orbits. 
Suppose that there exists a periodic orbit, say $P$. 
In this setting we concern the question to find more periodic orbits other than $P$ by using the notion of the braid type \cite{Boy}. 
Let $D_n$ be the $n$-punctured disk,  where $n$ is a period of $P$. 
Take an arbitrary orientation preserving homeomorphism $j: D \setminus P \rightarrow D_n$ and consider the composition 
$\hat{f}= j \circ f|_{D \setminus P} \circ j^{-1}: D_n \rightarrow D_n$. 
Let $[\hat{f}]$ be the isotopy class of $\hat{f}$.  
The {\it braid type of} $P$ {\it for} $f$, denoted by $bt(P,f)$, is the conjugacy class of $[\hat{f}]$ 
in the mapping class group $MCG(D_n)$ of $D_n$. 

By Nielsen--Thurston theory any homeomorphism of $D_n$ is isotopic to either periodic, reducible or pseudo--Anosov map. 
Since the statement is invariant under conjugacy, it makes sense to speak of the periodic, reducible, pseudo--Anosov braid type. 
The theory detects the complicated dynamics from the existence of periodic orbits of pseudo--Anosov braid type. 
For example if  $bt(P,f)$ is pseudo--Anosov there exists an infinitely many number of periodic orbits  with distinct periods for $f$. 
Moreover the logarithm of the dilatation of $bt(P,f)$ gives the lower bound 
of the topological entropy for $f$  \cite[Expos\'{e} 10]{FLP}. 
Recently it has been recognized that such complexity realizes global and efficient  particle mixings in fluid dynamics  \cite{BAS,KS}. 

Our interest is to show which braid types are forced by a given periodic orbit. 
For the study we use the language of the forcing relation on braid types. 
We denote the set of braid types of all periodic orbits for $f$ by $bt(f)$. 
Let $BT_n$ be the set of braid types of period $n$ orbits for all homeomorphisms of $D$, and 
$BT= \{ \beta \in BT_n\ |\ n \ge 1\ \}$.  
For an element $mc \in MCG(D_n)$, $[mc]$ denotes its conjugacy class. 
Following \cite{Han} we say that $f: D \rightarrow D$ {\it exhibits} $[mc] \in BT_n$ if 
there exists a periodic $n$ orbit for $f$ whose braid type is $[mc]$. 
We say that $[mc_1] \in BT$ {\it forces} $[mc_2] \in BT$, denoted by $[mc_1] \succeq [mc_2]$ provided that 
if a homeomorphism $f: D \rightarrow D$ exhibits $[mc_1]$,  then $f$ also exhibits $[mc_2]$. 
This relation $\succeq$ is a partial order on $BT$ \cite{Boy, Los}, and it is called the {\it forcing relation} or {\it forcing partial order}. 

Let $[mc_1] , [mc_2] \in BT$, and suppose that $[mc_1]$ is pseudo--Anosov. 
Results by Asimov--Franks \cite{AF} and Hall \cite{Hal1} give a strategy to determine whether $[mc_1]$ forces $[mc_2]$ or not. 
It holds that $[mc_1] \succeq [mc_2]$ if and only if the pseudo--Anosov map 
$\Phi_{mc_1} \in mc_1$ (which is thought as a homeomorphism of $D$) exhibits $[mc_2]$. 

For the study of braid types it is convenient to use geometric braids. 
There is a surjective homomorphism $\Gamma$ from the $n$--braid group $B_n$ to $MCG(D_n)$. 
We write $\sigma_i$, $i= 1, \cdots, n-1$ for the Artin generators of $B_n$. 
Any braid type is written by $[\Gamma(b)]$ for some braid $b \in B_n$. 
Simply we write $[b]$ for $[\Gamma(b)] \in BT_n$ when there is no confusion.  

\begin{figure}[htbp]
\begin{center}
\includegraphics[width=3.5in]{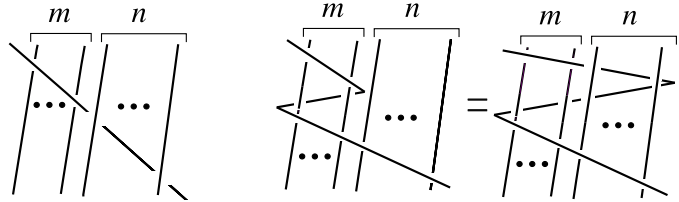}
\caption{Braids $\beta_{m,n}$ (left) and $\sigma_{m,n}$ (right). }
\label{fig_b-s_braid}
\end{center}
\end{figure}

This paper concerns the forcing partial order on the sets of braid types 
$\{[\beta_{m,n}]\}_{m,n \ge 1}$ and $\{[\sigma_{m,n}]\}_{m,n \ge 1}$ defined as follows. 
For any positive integers $m$ and $n$, let $\beta_{m,n}$ and $\sigma_{m,n}$ be 
the $m+n+1$--braids as in Figure~\ref{fig_b-s_braid}. 
The braid $\sigma_{m,n}$ can be written as $\sigma_{m,n}=\beta_{m,n} \xi$, where 
$\xi= \sigma_{m+n} \cdots \sigma_2 \sigma_1 \sigma_1 \sigma_2 \cdots \sigma_{m+n}$ (Figure~\ref{fig_b-s_braid}(right)).   
Each $\beta_{m,n}$ is pseudo--Anosov, and $\sigma_{m,n}$ is pseudo--Anosov if and only if $|m-n| \ge 2 $. 
These braids are concerned  in \cite{HK} from view point of  braids with small dilatation. 

The following is the main result of this paper: 
\begin{thm}
\label{thm_order1} 
For any $m,n \ge 1$ we have: 
\begin{description}
\item{(1)}
$[\beta_{m,n}] \succeq [\beta_{m+1,n}]$ and $[\beta_{m,n}] \succeq [\beta_{m,n+1}]$, 
\item{(2)} 
$[\beta_{m,n}] \succeq [\sigma_{m,\ell}]$ if $\ell \ge m+2$, and 
\item{(3)} 
$[\sigma_{m,n}] \succeq [\sigma_{m,\ell}]$ if $n \ge \ell \ge m+2$. 
\end{description}
\end{thm}
In particular fixing $m \ge 1$,  the relation $\succeq $ becomes a total order on each set of braid types $\{[\beta_{m,n}]\}_{n \ge 1}$, $\{[\beta_{n,m}]\}_{n \ge 1}$ and $\{[\sigma_{m,n}]\}_{n \ge m+2}$. 

\begin{figure}[htbp]
\begin{center}
\includegraphics[width=3in]{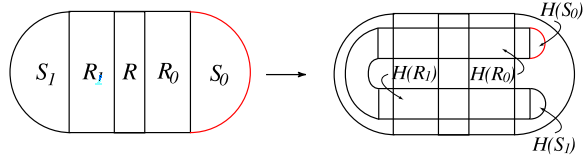}
\caption{Smale--horseshoe map $\mathtt{H}$.}
\label{fig_smale-h}
\end{center}
\end{figure}

{\it Horseshoe braid types} are those that can be realized by the periodic orbits 
for the Smale--horseshoe map $\mathtt{H}:D \rightarrow D$ (Figure~\ref{fig_smale-h}). 
This map is known to be a simple factor possessing  complicated dynamics. 
The following theorem says that $[\beta_{1,k}]$ forces any horseshoe braid type: 

\begin{thm} 
\label{thm_horseshoe} 
We have $[\beta_{1,k}] \succeq [mc]$ for any $k \ge 1$ and any horseshoe braid type $[mc]$. 
\end{thm}
Kolev shows that if $f$ has a period $3$ orbit $P$ whose braid type is pseudo--Anosov, $f$ has a periodic orbit of every period \cite{Kol}. 
This is a best possible analogous result of the Li--York theorem \cite{LY} (or the special case of the Sharkovskii theorem), although  
the theorem does not say which braid types can be realized by the period $3$ orbit. 
Note that the braid type for a fixed point or a period $2$ orbit is unique. 
A question is which braid type for a period $n$ orbit $(n \ge 4)$ is forced by a period $3$ orbit of pseudo--Anosov type. 
Theorems~\ref{thm_order1} and \ref{thm_horseshoe} together with a Handel's result \cite{Han} gives an answer: 

\begin{cor}
\label{cor_period3}
Let $f:D \rightarrow D$ be an orientation preserving homeomorphism. 
Suppose that $f$ has a period $3$ orbit whose braid type is pseudo--Anosov. 
Then we have: 
\begin{description}
\item{(1)} 
$bt(f) \supset \{[\beta_{m,n}]\ |\ m,n \ge 1\}$, and 
\item{(2)} 
 $bt(f) \supset bt(\mathtt{H}) \supset \{[\sigma_{m,n}]\ |\ n \ge m+2\}$. 
\end{description}
\end{cor}

{\sc Acknowledgements.} 
The author thanks Yoshihiro Yamaguchi and Kiyotaka Tanikawa for giving her a motivation to  study families of braids in this paper. 
The author also thanks Takashi Matsukoka for a great deal of encouragement.
The author is partially supported by Grant--in--Aid for Young Scientists (B) (No. 17740094), 
The Ministry of Education, Culture, Sports, Science and Technology, Japan.

\section{Preliminaries} 
\label{section_preliminaries}

In Section~\ref{subsection_pA-dil-for} we review the classification theorem of surface homeomorphisms by Nielsen--Thurston 
and a result on the relation between the forcing partial order and the dilatation of braids by Los. 
Section~\ref{subsection_fibered} introduces a fibered surface for a given graph, and 
it describes a criterion for determining whether  a braid $b$ is pseudo--Anosov or not by Bestvina--Handel. 
Under the assumption that $\mathfrak{g}$ is the induced graph map for a braid $b$ satisfying the Bestvina-Handel's condition,  
in Section~\ref{subsection_thick-pA} we define the reduced graph map $\mathfrak{g}^{\mathrm{red}}$, and 
we give a relation between periodic orbits for the thick graph map associated to 
$\mathfrak{g}^{\mathrm{red}}$ and those for the pseudo--Anosov map $\Phi_b \in \Gamma(b)$. 
Section~\ref{subsection_Smale} explains the dynamics on the Smale--horseshoe map 
can be described by the shift map on the symbol space, and it shows that 
the braids $\sigma_{m,n}$ ($n \ge m+2$) arise as braid types of periodic orbits.

\subsection{Pseudo-Anosov braids}
\label{subsection_pA-dil-for}

We introduces three kinds of  homeomorphisms. 
A homeomorphism $\Phi: D_n \rightarrow D_n$ is: 
\begin{description}
\item{}
{\it periodic} if some power of $\Phi$ is the identity map. 
\item{}
{\it reducible} if there is a $\Phi$--invariant closed $1$-submanifold whose complementary components in 
$D_n$ have negative Euler characteristic. 
\item{}
{\it pseudo--Anosov} if there is a constant $\lambda >1$ and a pair of transverse measured foliations 
$(\mathcal{F}^s, \mu^s)$ and   $(\mathcal{F}^u, \mu^u)$ 
such that $\Phi(\mathcal{F}^s, \mu^s) =(\mathcal{F}^s, \lambda^{-1} \mu^s )$ and 
$\Phi(\mathcal{F}^u, \mu^u) =(\mathcal{F}^u, \lambda \mu^u)$. 
\end{description}
$\mathcal{F}^{s}$ and $\mathcal{F}^{u}$ are called the
 {\it stable} and {\it unstable foliation} or the {\it invariant foliations}. 
They have a finitely many number of singularities, and 
the set of singularities of $\mathcal{F}^s$ equals that of $\mathcal{F}^u$. 
The number $\lambda= \lambda(\Phi)>1$ is called the {\it dilatation} for $\Phi$. 

We say that $mc \in MCG(D_n)$ is {\it periodic} ({\it reducible}, {\it pseudo--Anosov} resp.) 
if it contains a periodic map (reducible map, pseudo--Anosov map resp.) as a representative homeomorphism. 
An element $mc \in MCG(D_n)$ is called {\it irreducible} if it is not reducible. 

\begin{thm} \cite{FLP}
Any irreducible element $mc \in MCG(D_n)$ is periodic or pseudo--Anosov. 
If $mc$ is pseudo--Anosov, then the pseudo--Anosov map of $mc$ is unique up to conjugacy. 
\end{thm}
The Nielsen--Thurston type (i.e, periodic, reducible, pseudo--Anosov) for $mc$ 
is invariant under conjugacy. When $mc$ is pseudo--Anosov,  
the {\it dilatation} $\lambda(mc)$ for  $mc$ is defined by $\lambda(\Phi_{mc} )$ 
for the pseudo--Anosov map $\Phi_{mc} \in mc$. 
This number is also invariant under conjugacy.

\begin{figure}[htbp]
\begin{center}
\includegraphics[width=2.3in]{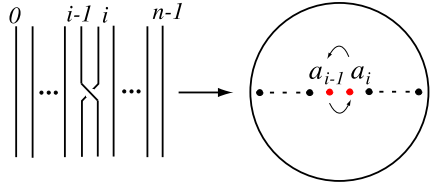}
\caption{$\Gamma: B_n \rightarrow MCG(D_n)$. }
\label{fig_homo}
\end{center}
\end{figure}

Let $A_n= \{a_0, \cdots, a_{n-1}\}$ be a set of $n$--points in the interior of $D$. 
Suppose that $a_0, \cdots, a_{n-1}$ lie on the horizontal line through the center of the disk from 
the left to the right, and put $D_n = D \setminus A_n$. 
Let $D_i$, $i = 0, \cdots, n-1$ be the closed disk 
which contains $a_{i-1}$ and $a_{i}$ and no other points of  $A_n$. 
We define a  homomorphism 
$\Gamma: B_n \rightarrow MCG(D_n)$ as follows: 
For the Artin generators $\sigma_i$, $i=1, \cdots, n-1$, 
$\Gamma(\sigma_i)$ is represented by a homeomorphism of $D_n$ 
which fixes the exterior of $D_i$ and rotates in the inside of $D_i$ by $180$ degrees in the counter--clockwise direction 
so that $a_{i-1}$ is interchanged with $a_i$ (Figure~\ref{fig_homo}). 
The kernel of $\Gamma$ is the center of $B_n$ which is generated by a full twist braid 
$(\sigma_1 \sigma_2 \cdots \sigma_{n-1})^n$ \cite{Bir}. 
We say that a braid $b \in B_n$ is {\it pseudo--Anosov} ({\it periodic}, {\it reducible} resp.) 
if $\Gamma(b) \in MCG(D_n)$ is pseudo--Anosov (periodic, reducible resp.). 
We define the {\it dilatation} $\lambda(b)$ for the pseudo--Anosov braid $b$ by $\lambda(\Gamma(b))$. 

Recall that we  write $[b]$ for $[\Gamma(b)]$. 
One of the relation between the forcing partial order and the dilatation is as follows: 

\begin{thm}\cite{Los}
\label{thm_Los}
Suppose that $b_1$ and $b_2$ are pseudo--Anosov braids. 
If $[b_1] \succeq [b_2]$ with $[b_1] \ne [b_2]$, then $\lambda(b_1) > \lambda(b_2)$. 
\end{thm}

\subsection{Graphs, fibered surfaces and graph maps} 
\label{subsection_fibered}

Let  $G$ be a finite graph embedded on an orientable surface $F$. 
In this paper, we assume that an edge of $G$ is closed, and 
let $\mathcal{E}_{\mathrm{ori}}(G)$ be the set of oriented edges of $G$, 
$\mathcal{E}(G)$ the set of unoriented edges, and $\mathcal{V}(G)$ the set of vertices.  
We denote the oriented edge with the initial vertex $v_I$ and the terminal vertex $v_T$ by $e(v_I, v_T)$. 
Let $\overline{e}$ be the same edge as $e$ with opposite orientation. 
A continuous map $\mathfrak{g}: G \rightarrow G$ is called a {\it graph map}. 

\begin{figure}[htbp]
\begin{center}
\includegraphics[width=2.6in]{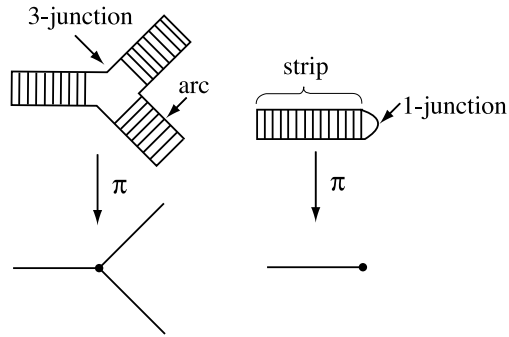}
\caption{Fibered surface.}
\label{fig_fiber-surface}
\end{center}
\end{figure}

One can associate a {\it fibered surface} ${\Bbb F}(G) \subset F$ 
with a projection $\pi: {\Bbb F} (G) \rightarrow G$ (Figure~\ref{fig_fiber-surface}). 
The fibered surface ${\Bbb F}(G)$ is decomposed into arcs and into polygons modelled on $k$--junctions, $k \ge 1$. 
The arcs and the $k$--junctions are called  {\it decomposition elements}. 
Under $\pi$, the preimage of each vertices of valence $k$ is the $k$--junction, and the closure of the preimage of each open edge is the strip (fibered by arcs) 
which is the closure of  the one of the complementary components of the union of all junctions.

Recall that $A_n= \{a_0, \cdots, a_{n-1}\}$ is a set of $n$--points in  $D$ and $D_n = D \setminus A_n$. 
In this section we have no assumption on the location of $A_n$. 
Let $P_i$ be  a small circle centered at $a_i$ such that 
no other points of $A_n$ is contained in the disk bounded by $P_i$. 
We set $P= \displaystyle \bigcup_{i=0}^{n-1} P_i$. 
Choose a finite graph $G$ embedded on $D_n$ that is homotopy equivalent to $D_n$ 
such that $P$ is a subgraph of $G$ and $G$ has no vertices of valence $1$ or $2$. 

Let $f: D_n \rightarrow D_n$ be a representative homeomorphism of $\Gamma(b) \in MCG(D_n)$. 
A fibered surface ${\Bbb F}(G)$ {\it carries} $f$ if $f$ maps each decomposition element of ${\Bbb F}(G)$  into a decomposition element and 
each junction into a junction. 
When ${\Bbb F}(G)$ carries $f$, $f$ induces a  piecewise linear graph map 
$\mathfrak{g}=\mathfrak{g}_f: G \rightarrow G$ which represents the correspondence of 
vertices and edge paths determined by $f$. 
(Thus $\mathfrak{g}$ sends vertices to vertices and each edge to an edge path.) 
We may assume that $P$ is invariant under $\mathfrak{g}$ without loss of generality. 

Suppose that a  fibered surface ${\Bbb F}(G)$ carries $f: D_n \rightarrow D_n$  
of $\Gamma(b)$ with the induced graph map $\mathfrak{g}: G \rightarrow G$. 
In this case we say that  $\mathfrak{g}$ is the {\it induced graph map for }$b$. 
Let pre$P$ be the set of edges $e \in \mathcal{E}(G)$ such that 
$\mathfrak{g}^k(e)$ is contained in $P$ for some $k \ge 1$. 
The graph map $\mathfrak{g}$ defines a  square and non negative integral matrix called 
the {\it transition matrix} ${\mathcal T}_{\mathfrak{g}}^{\mathrm{tot}}$ 
whose $(i,j)$th entry is given by the number of times 
 that the image of the $j$th edge of $\mathcal{E}(G)$ under $\mathfrak{g}$ 
passes through the $i$th edge of $\mathcal{E}(G)$. 
Then ${\mathcal T}_{\mathfrak{g}}^{\mathrm{tot}}$ is of the form 
$${\mathcal T}_{\mathfrak{g}}^{\mathrm{tot}}= 
\left(\begin{array}{ccc}
{\mathcal P} & {\mathcal A} & {\mathcal B} 
\\
0 & {\mathcal Z} & {\mathcal C} \\0 & 0 & {\mathcal T}
\end{array}\right),$$
where ${\mathcal P}$ and $ {\mathcal Z} $ are the transition matrices associated to 
$P$ and pre$P$ respectively, and 
${\mathcal T}$ is the transition matrix associated to the rest of edges called {\it real edges}. 
The matrix ${\mathcal T}$ is called the {\it transition matrix with respect to the real edges}. 
The spectral radius of ${\mathcal T}$ is denoted by $\lambda({\mathcal T})$. 

%

A graph map $\mathfrak{g} : G \rightarrow G$ is {\it efficient} if for any $e \in  \mathcal{E}_{ori}(G)$ and any $k \ge  0$,
$\mathfrak{g}^k(e) = e_{k,1} e_{k,2} \cdots e_{k,j}$ satisfies $\overline{e_{k,i}} \neq e_{k,i+1}$ for all $i=1,\dots,j-1$.

A non negative square matrix  $M$ is {\it irreducible} if for every set of indices $i,j$, there is a positive integer  $n_{i,j} $ such that the $(i,j)$th entry of 
$M^{n_{i,j}}$ is strictly positive.  

\begin{thm} 
\cite{BH} 
\label{BH-thm}
Let $b \in B_n$ and let $\mathfrak{g}:G \rightarrow G$ the induced graph map for $b$.   
Suppose that 
\begin{description}
\item{(BH:1)} $\mathfrak{g}$ is efficient, and 
\item{(BH:2)} the transition matrix $\mathcal{T}$ with respect to the real edges is 
irreducible with $\lambda(\mathcal{T})>1$.
\end{description}
Then $b$ is pseudo--Anosov with dilatation $\lambda(\mathcal{T})$.
\end{thm}
An idea of the proof is as follows. 
The {\it train track} $\tau \subset D_n$ associated to $\mathfrak{g}$ is obtained by the \lq\lq smoothing" of $G$, and as a result 
the {\it train track map} $\mathfrak{g}_{\tau}: \tau \rightarrow \tau$ can be defined. 
If $\mathfrak{g}: G \rightarrow G$ satisfies (BH:1) and (BH:2), 
one can construct the pseudo--Anosov map $\Phi_{b} \in \Gamma(b)$ explicitly by using $\mathfrak{g}_{\tau}$,  and hence $b $ is pseudo--Anosov. 
For more details see \cite[Section 3.3]{BH}

\subsection{Thick graph maps and pseudo--Anosov maps} 
\label{subsection_thick-pA}.

\begin{figure}[htbp]
\begin{center}
\includegraphics[width=3in]{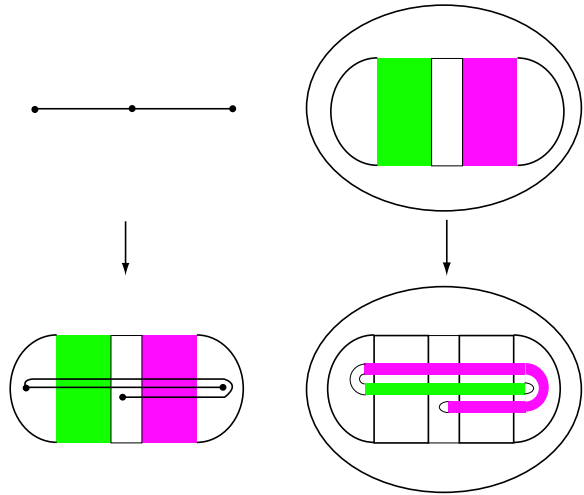}
\caption{Embedding (left) and  its  thick graph map (right).}
\label{fig_per3_horseshoe}
\end{center}
\end{figure}

Let $T$ be a finite tree embedded on $D$ and $\widetilde{\mathfrak{g}}: T \rightarrow  {\Bbb F}(T)$ an embedding such that 
it maps a vertex to a junction and the image of each edge is transverse to arcs of ${\Bbb F}(T)$. 
A homeomorphism $g: D (\supset {\Bbb F}(T)) \rightarrow D$ is a {\it thick graph map associated to} $\widetilde{\mathfrak{g}}$ 
if $g$ satisfies the following conditions: 
\begin{itemize}
\item 
$g$ maps each decomposition element of ${\Bbb F}(T)$ into a decomposition element and each junction into a junction. 
\item 
$g$ contracts the vertical direction of each strip of ${\Bbb F}(T)$ uniformly and expands the horizontal direction of each strip uniformly. 
\item 
$g({\Bbb F}(T))$ is a fibered surface of the tree $\widetilde{\mathfrak{g}}(T)$. 
\end{itemize}
For example see Figure~\ref{fig_per3_horseshoe}. 
Although the thick graph map $g$ is not unique, it is determined uniquely 
on  the invariant set $\Lambda= \displaystyle \bigcap_{j \in {\bf Z}} g^j({\Bbb F}(T))$ 
under $g$ in a sense of the symbolic dynamics \cite{ALM}. 

\begin{figure}[htbp]
\begin{center}
\includegraphics[width=6in]{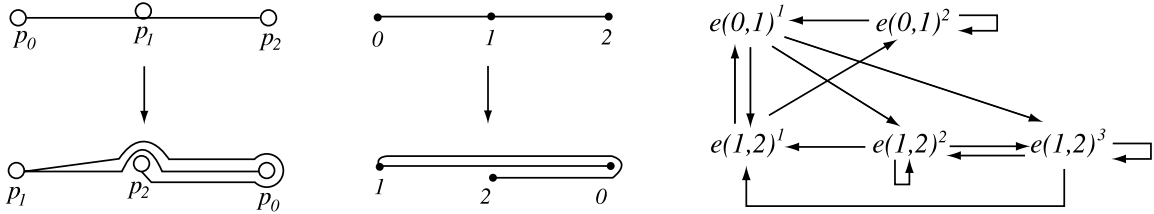}
\caption{Induced graph map for $\sigma_1 \sigma_2^{-1}$ (left), 
its reduced graph map (center) and transition graph (right). }
\label{fig_graph_map11}
\end{center}
\end{figure}

In the rest of this section, let $\mathfrak{g}: G \rightarrow G$ be the induced graph map for $b \in B_n$  satisfying (BH:1) and (BH:2). 
The {\it reduced graph} $G^{\mathrm{red}}$ for $G$ is the tree 
obtained from $G$ by collapsing each peripheral edge $P_i$ of $P$ to a vertex $v_i$ labeled $i$.  
Since $P$ is invariant under $\mathfrak{g}$, 
a piecewise linear graph map $\mathfrak{g}^{\mathrm{red}}: G^{\mathrm{red}} \rightarrow G^{\mathrm{red}}$, 
called the {\it reduced graph map for} $\mathfrak{g}$,  can be defined such that 
$\mathfrak{g}^{\mathrm{red}}$ represents the correspondence of vertices and edge paths determined by $\mathfrak{g}$, 
see Figure ~\ref{fig_graph_map11}(left/center). 
Because $\mathfrak{g}$ is the induced graph map for $b$, ${\Bbb F}(G)$ carries some homeomorphism $f$ of $\Gamma(b)$. 
Thus $\mathfrak{g}^{\mathrm{red}}: G^{\mathrm{red}} \rightarrow G^{\mathrm{red}}$ recovers the embedding 
$\widetilde{\mathfrak{g}^{\mathrm{red}}}: G^{\mathrm{red}} \rightarrow {\Bbb F}(G^{\mathrm{red}})$. 
 The thick graph map $g: D \rightarrow D$ {\it associated to} $\mathfrak{g}^{\mathrm{red}}$ means that the one  associated to $\widetilde{\mathfrak{g}^{\mathrm{red}}}$.

Assume that the closed braid of $b$ is a knot. 
Then  the thick graph map $g: D \rightarrow D$ associated to $\mathfrak{g}^{\mathrm{red}}$ has a period $n$ orbit, say $P_n= \{p_0, \cdots, p_{n-1}\}$ such that a point $p_i$ is in a junction $\pi^{-1}(v_i)$. 
Notice that   the braid type of $P_n$ for $g$ is $[b]$. 
We call $P_n$ the {\it representative orbit} for $g$ (associated to $\mathfrak{g}^{\mathrm{red}}$). 

Results by Asimov--Franks and Hall tell us that 
$$bt(\Phi_b) =\{[b'] \in BT\ |\ [b] \succeq [b']\}.$$
Hence we have $bt(g) \supset  bt(\Phi_{b})$ since $[b] \in bt(g)$. 
To give elements of $bt(g)$ which belong to $bt(\Phi_{b})$, we introduce terminology. 
A periodic point $x \in G^{\mathrm{red}}$ for $\mathfrak{g}^{\mathrm{red}}$ is called {\it regular} 
if $x \notin \mathcal{V}(G^{\mathrm{red}})$. 
Since $\mathfrak{g}^{\mathrm{red}}(\mathcal{V}(G^{\mathrm{red}})) \subset \mathcal{V}(G^{\mathrm{red}})$, 
if $x$ is regular, $(\mathfrak{g}^{\mathrm{red}})^i (x)$ is also regular  for each $i \ge 0$.  
Hence it makes sense to speak of the regular periodic orbit. 
Since the number of the element of $\mathcal{V}(G^{\mathrm{red}})$ is finite, 
it is easy to check that a given periodic orbit is regular 
or not. 

Next we define a subdivision $G^{\mathrm{red}}_s$ of $G^{\mathrm{red}}$ as follows: 
Let $e$ be an edge of $ \mathcal{E}(G^{\mathrm{red}})$ such that 
the edge path $\mathfrak{g}^{\mathrm{red}}(e)$ is given by 
$ f_1 f_2 \cdots f_k$, $f_i \in \mathcal{E}(G^{\mathrm{red}})$. 
Subdivide $e$ into subedges $e^1, e^2, \cdots, e^k$ 
so that $\mathfrak{g}^{\mathrm{red}}(e^i)= f_i$, $i = 1, \cdots, k$. 
(Thus $\mathfrak{g}^{\mathrm{red}}(e^1 e^2 \cdots e^k)= f_1 f_2 \cdots f_k$ as an edge path.)
Let $E_1, \cdots, E_{\ell} \in \mathcal{E}(G^{\mathrm{red}}_s)$ be all edges of $G^{\mathrm{red}}_s$. 
The {\it transition graph} $\Xi_{\mathfrak{g}}$ is the oriented graph with vertices 
$E_1, \cdots, E_{\ell}$ and an oriented edge from $E_i$ to $E_j$ 
if $\mathfrak{g}^{\mathrm{red}}(E_i)$ passes through $E_j$. 
Note that from the definition of the subdivision, 
$\mathfrak{g}^{\mathrm{red}}(E_i)$ passes through $E_j$ at most one times. 
For example consider  the reduced graph map given in 
Figure~\ref{fig_graph_map11}(center), and in this case $e(0,1)$ 
is subdivided into $e(0,1)^1$ and $e(0,1)^2$, and 
$e(1,2)$ is subdivided into $e(1,2)^1$, $e(1,2)^2$ and $e(1,2)^3$. 
Since 
\begin{eqnarray*}
\mathfrak{g}^{\mathrm{red}}(e(0,1))&=& \mathfrak{g}^{\mathrm{red}}(e(0,1)^1 e(0,2)^2)= 
e(2,1)e(1,0) \hspace{2mm}\mbox{and} 
\\
\mathfrak{g}^{\mathrm{red}}(e(1,2))&=& \mathfrak{g}^{\mathrm{red}}(e(1,2)^1 e(1,2)^2 e(1,2)^3)
= e(0,1)e(1,2)e(2,1), 
\end{eqnarray*}
we have the transition graph shown in Figure~\ref{fig_graph_map11}(right). 

Each closed path of $\Xi_{\mathfrak{g}}$ gives a periodic orbit for 
the thick graph map $g$ associated to $\mathfrak{g}^{\mathrm{red}}$: 

\begin{lem}
\label{lem_symbol}
Let  $E_0 \rightarrow \cdots \rightarrow E_{s-1} \rightarrow E_0$, $E_i \in \mathcal{V}(\Xi_{\mathfrak{g}})$ 
be a closed path of $\Xi_{\mathfrak{g}}$. 
Then 
\begin{description}
\item{(1)}  
there exists a periodic point $x_0 \in E_0$ for 
$\mathfrak{g}^{\mathrm{red}}$ 
such that  $(\mathfrak{g}^{\mathrm{red}})^s(x_0)= x_0$ and 
$x_i= (\mathfrak{g}^{\mathrm{red}})^i(x_0) \in E_i$ 
for each $i \in \{0, \cdots, s-1\}$, and 
\item{(2)}
 there exists a periodic point $\widehat{x_0} \in \pi^{-1}(E_0)$ for $g$ 
 associated to $\mathfrak{g}^{\mathrm{red}}$ 
such that  $g^s(\widehat{x_0})= \widehat{x_0}$ and 
$\widehat{x_i}=g^{i}(\widehat{x_0}) \in \pi^{-1}(E_i)$ for each $i \in \{0, \cdots, s-1\}$.
\end{description}
\end{lem}

Proof. 
This can be shown by the symbolic dynamics (for example see \cite{ALM}). $\Box$
\medskip

The word $E_0 E_1 \cdots E_{s-1}$ ($\pi^{-1}(E_0) \cdots \pi^{-1}(E_{s-1})$ resp.) in Lemma~\ref{lem_symbol} 
is said to be the {\it itinerary} of $x_0$ ($\widehat{x_0}$ resp.). 
\medskip


The symbol $\mathcal{O}_{f}(x)$ denotes the periodic orbit for a periodic point $x$ for a map $f$. 
Since $\Phi_{b} $ (and also Markov partition of $\Phi_b$)
is constructed via $\mathfrak{g}: G \rightarrow G$, 
there is a natural correspondence between periodic orbits for $\Phi_b$ and those for $g$. 
If $x_0$ is the periodic point for $\mathfrak{g}^{\mathrm{red}}$ associated with a closed path 
$E_0 \rightarrow \cdots \rightarrow E_{s-1} \rightarrow E_0 $ of $\Xi_{\mathfrak{g}}$ 
in the sense of Lemma~\ref{lem_symbol},  there is a periodic point $\widetilde{x_0}$ for $\Phi_{b}$ such that 
$(\Phi_{b})^s(\widetilde{x_0})=\widetilde{x_0}$ and 
$\widetilde{x_i}=(\Phi_{b})^i(\widetilde{x_0})$ is in a {\it Markov box} labeled $E_i$ for each 
$i \in \{0, \cdots, s-1\}$. 
If $x_0$ is regular, then the periodic orbit of $\widetilde{x_0}$ lie  on the regular (non--singular) leaves of the 
stable and unstable foliations for $\Phi_{b}$. 
Then the construction of two maps $\Phi_b$ and $g$ implies that 
$bt(\mathcal{O}_{\Phi_b}(\widetilde{x_0}), \Phi_b)= bt(\mathcal{O}_{g} (\widehat{x_0}), g)$. 
(In particular, the period of $\widetilde{x_0}$ for $\Phi_{b}$ 
equals that of $\widehat{x_0}$ for $g$.) 
For more details see \cite[Section~3.3]{BH}. 
Thus we have: 
\begin{lem}
\label{lem_symbol2} 
In Lemma~\ref{lem_symbol} if $x_0$ is regular, then 
$bt(\Phi_b) \ni bt(\mathcal{O}_{g} (\widehat{x_0}), g) $, and hence 
$$ [b] \succeq bt(\mathcal{O}_{g} (\widehat{x_0}), g) .$$  
\end{lem}

If $x_0$ is not regular, $g$ has a periodic point $x_0'$ in the junction $\pi^{-1}(x_0)$. 
Then $bt(\mathcal{O}_g(x_0'), g)=bt(\mathcal{O}_{\Phi_b}(\widetilde{x_0}), \Phi_b)$ 
from the construction of two maps.  
However it is not true in general that the period of $x_0'$ for $g$ equals that of $\widehat{x_0}$ for $g$. 
In this case $bt(\mathcal{O}_{\Phi_b}(\widetilde{x_0}), \Phi_b) \ne bt(\mathcal{O}_{g} (\widehat{x_0}), g)$.  

\subsection{Smale--horseshoe map}
\label{subsection_Smale}

The Smale--horseshoe map $\mathtt{H}: D \rightarrow D$ is a diffeomorphism such that 
the action of $\mathtt{H}$ on three rectangles $R_0,R_1$ and $R$ and two half disks 
$S_0,S_1$ is given in Figure~\ref{fig_smale-h}. 
The restriction $\mathtt{H}|_{R_i}$, $i=0,1$ is an affine map such that 
$\mathtt{H}$ contracts $R_i$ vertically and stretches horizontally, and 
 $\mathtt{H}|_{S_i}$, $i=0,1$ is a contraction map. 

The set  $\Omega = \displaystyle\bigcap_{j \in {\bf Z}} \mathtt{H}^j (R_0 \cup R_1)$ is invariant under $\mathtt{H}$, and 
$\mathtt{H}|_{\Omega}: \Omega \rightarrow \Omega$ is conjugate to the {\it shift map}  
$\sigma: \Sigma_2=\{0,1\}^{\Bbb Z} \rightarrow \Sigma_2$, where 
\begin{eqnarray*}
\sigma(*** w_{-1}\cdot w_0 w_1 ***) = (*** w_{-1}w_0 \cdot w_1 ***),
\quad w_j \in \{0,1\}. 
\end{eqnarray*}
The conjugacy $\mathcal{K}: \Omega \rightarrow \Sigma_2$ is given by 
\begin{eqnarray*}
\mathcal{K}(x) &=& (\cdots \mathcal{K}_{-1}(x) \mathcal{K}_0(x) \mathcal{K}_1(x) \cdots),
\hspace{2mm}
\mbox{where}
\end{eqnarray*}
\[
\mathcal{K}_j(x) =
\left\{
\begin{array}{ll}
0 \hspace{3mm}\  \mbox{if\ }&\mathtt{H}^j(x) \in R_0,
\\
1 \hspace{3mm}\  \mbox{if\ }&\mathtt{H}^j(x) \in R_1.
\end{array}
\right.
\]
If $x$ is a period $k$ point, the word $ \mathcal{K}_0(x) \mathcal{K}_1(x) \cdots, \mathcal{K}_{k-1}(x)$ is called the {\it code} for $x$. 
Modulo cyclic permutation,  $\mathcal{K}_0(x) \mathcal{K}_1(x) \cdots, \mathcal{K}_{k-1}(x)$ 
is said to be the {\it code} for the periodic orbit $\mathcal{O}_{\mathtt{H}}(x)$. 
We say that $[b]$, $b \in B_n $   is a {\it horseshoe braid type} 
if there is a period $n$ orbit for $\mathtt{H}$ whose braid type is $[b]$. 
For the study of the forcing partial order on the set of horseshoe braid types, 
see \cite{dCH,Hal2}. 

\begin{figure}[htbp]
\begin{center}
\includegraphics[width=2.4in]{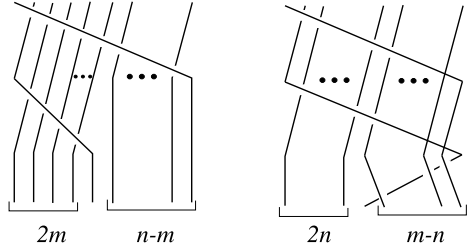}
\caption{Braid $\sigma'_{m,n}$: case $n \ge m$ (left), case $m \ge n$ (right).}
\label{fig_s-braid}
\end{center}
\end{figure}

The argument in \cite[Section~3.2]{HK} shows that 
$\sigma_{m,n}$ is conjugate to $\sigma'_{m,n}$ given in Figure~\ref{fig_s-braid}. 
It is not hard to see that when $n \ge m+2$, 
$[\sigma'_{m,n}](=[\sigma_{m,n}])$ is a horseshoe braid type 
such that the corresponding periodic orbit for $\mathtt{H}$ has a code 
$1 \underbrace{0 \cdots 0}_{n-1} 1 \underbrace{0 \cdots 0}_{m}$ or
$1 \underbrace{0 \cdots 0}_{n-1} 1 \underbrace{0 \cdots 0}_{m-1}1$. 
For example $[\sigma_{1,3}']= [\sigma_1 \sigma_2 \sigma_3 \sigma_4 \sigma_1 \sigma_2]$ 
is the horseshoe braid with the code $10010$ or $10011$.

\section{Proof of theorems} 
\label{section_proofs}

\begin{figure}[htbp]
\begin{center}
\includegraphics[width=4in]{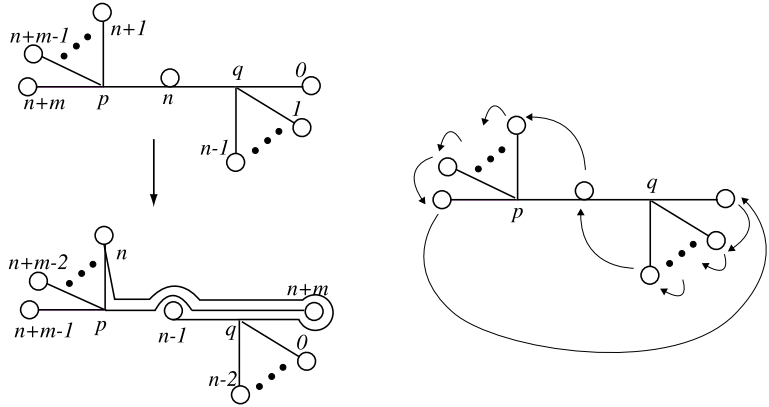}
\caption{$\mathfrak{g}_{m,n}: G_{m,n} \rightarrow G_{m,n}$ (left), 
transition of peripheral edges (right).}
\label{fig_b-graph_map}
\end{center}
\end{figure}

Let $\mathfrak{g}_{m,n}: G_{m,n} \rightarrow G_{m,n}$ be the graph map as in 
Figure~\ref{fig_b-graph_map}(left). 
 We label the vertices of $G_{m,n}$ which lie on the peripheral edges, 
$0,1, \cdots, n+m$ from the right to the left. 
Other vertices $p$ and $q$ of $G_{m,n}$ have valences  $m+1$ and $n+1$ respectively. 
 This is the induced graph map for $\beta_{m,n}$ satisfying (BH:1) and (BH:2) (\cite{HK}). 
 Hence $\beta_{m,n}$ is pseudo--Anosov for all $m,n \ge 1$. 
 Since $\beta_{m,n}^{-1}$ is conjugate to $\beta_{n,m}$, 
we have $\lambda(\beta_{m,n})= \lambda(\beta_{n,m})$. 
Figure~\ref{fig_b-graph_map}(right) indicates the transition of peripheral  edges. 

\begin{figure}[htbp]
\begin{center}
\includegraphics[width=4.5in]{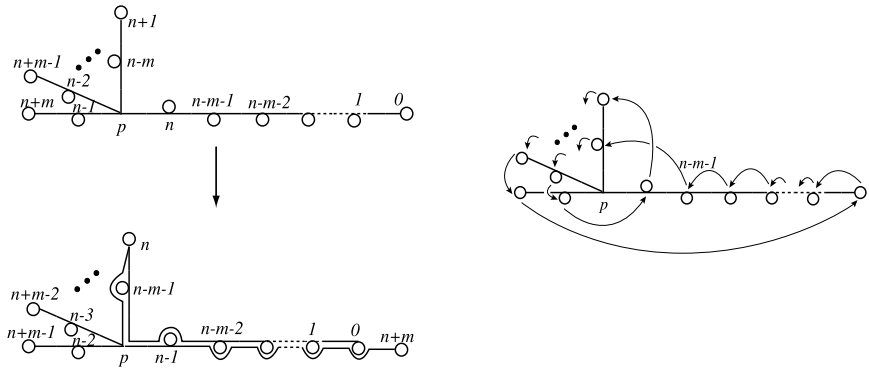}
\caption{$\mathfrak{h}_{m,n}: H_{m,n} \rightarrow H_{m,n}$ (left), transition of 
peripheral edges (right). }
\label{fig_s-graph_map}
\end{center}
\end{figure}

Now we turn to $\sigma_{m,n}$. 
For $n \ge m+2$, let $\mathfrak{h}_{m,n}: H_{m,n} \rightarrow H_{m,n}$ be the graph map as in Figure~\ref{fig_s-graph_map}(left). 
This is the induced graph map for $\sigma_{m,n}'$ 
in Figure~\ref{fig_s-braid}(left) satisfying (BH:1) and (BH:2) (\cite{HK}). 
Hence $\sigma_{m,n}$ is pseudo--Anosov in this case. 

\medskip

Proof of Theorem~\ref{thm_order1}(1). 
It is enough to show that $[\beta_{m,n}] \succeq [\beta_{m+1,n}]$, for 
$\beta_{m,n}$ is conjugate to $\beta_{n,m}^{-1}$ and in general 
$[b] \succeq [c]$ if and only if $[b^{-1}] \succeq [c^{-1}]$. 

We consider the induced graph map  $\mathfrak{g}_{m,n}$ for $\beta_{m,n}$. 
The transition graph $\Xi_{\mathfrak{g}_{m,n}}$ has vertices 
\begin{eqnarray*}
e(q,0),  e(q,1), \cdots, e(q,n-1), 
e(q,n)^1 ,  e(q,n)^2, e(q,n)^3, e(q,n)^4,  e(q,n)^5, 
\\
e(p,n) ,  e(p,n+1), \cdots, e(p,n+m-1), 
e(p,n+m)^1, e(p,n+m)^2, e(p,n+m)^3. 
\end{eqnarray*}
Let $g_{m,n}$ be the thick graph map associated to $\mathfrak{g}_{m,n}^{\mathrm{red}}$, 
$P_{m,n}$ the representative orbit for $g_{m,n}$. 
Note that $bt(P_{m,n}, g_{m,n})= [\beta_{m,n}]$. 
We consider the closed path $\mathcal{C}$ of $\Xi_{\mathfrak{g}_{m,n}}$ of  length $m+n+2$ 
such that 
\begin{eqnarray*}
E_0 &=& e(q,0) \rightarrow E_1 = e(q,1) \rightarrow 
\cdots \rightarrow E_{n-1}= e(q, n-1) \rightarrow E_n= e(q,n)^4 \rightarrow 
\\
E_{n+1}&=& e(p,n) \rightarrow E_{n+2}= e(p,n+1) \rightarrow \cdots 
\rightarrow E_{n+i}= e(p,n+i-1) \rightarrow 
\\
&\cdots& \rightarrow
E_{n+m}= e(p,n+m-1) \rightarrow E_{m+n+1}= e(p,n+m)^3 \rightarrow E_0 .
\end{eqnarray*}
Take a periodic point $x_0 \in e(q,0)$ for $\mathfrak{g}^{\mathrm{red}}_{m,n}$ 
given in Lemma~\ref{lem_symbol}. 
Then $x_n = (\mathfrak{g}_{m,n}^{\mathrm{red}})^n(x_0) \in E_n=e(q,n)^4$. 
Since $e(q,n)^4$ is a proper subedge of $e(q,n)$, $x_n $ is regular. 
Now we claim that the period of the orbit of $x_0$ is $m+n+2$.
Because $E_i \ne E_j$ ($i \ne j$) in $\mathcal{C}$, $\mathcal{C}$ is not a repetition of some closed subpath. 
Since $x_n$ is regular, $x_i \in E_i$ ($0 \le i \le m+n+1$) does not lie on the boundary of $E_i$. 
This implies that  the length of $\mathcal{C}$ equals the  period of $x_0$. 

\begin{figure}[htbp]
\begin{center}
\includegraphics[width=5.5in]{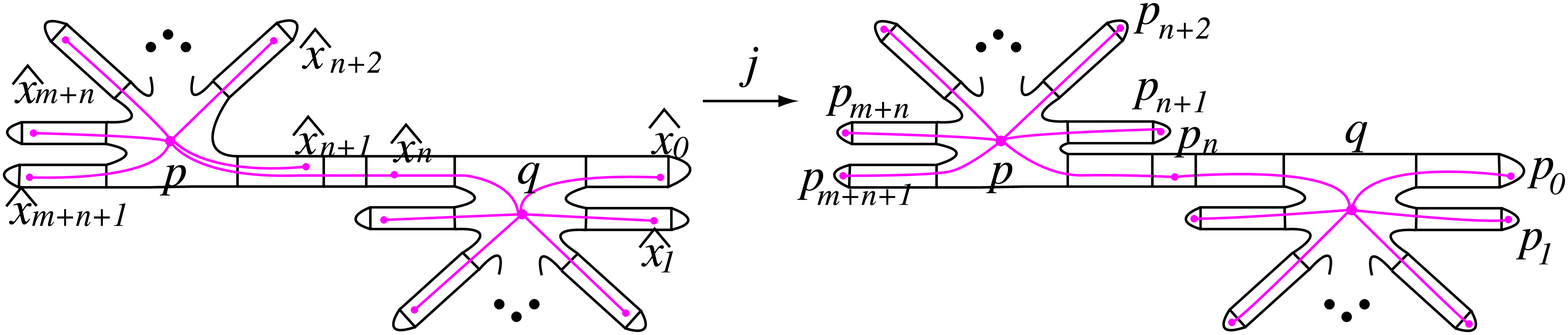}
\caption{$j$ sending $\widehat{G}$ 
to $G_{m+1,n}^{\mathrm{red}}$: 
$\widehat{G} \subset {\Bbb F}(G_{m,n}^{\mathrm{red}})$ (left), $G_{m+1,n}^{\mathrm{red}} \subset {\Bbb F}(G_{m+1,n}^{\mathrm{red}})$ (right).}
\label{fig_proof_isotopy}
\end{center}
\end{figure}

\begin{figure}[htbp]
\begin{center}
\includegraphics[width=4.8in]{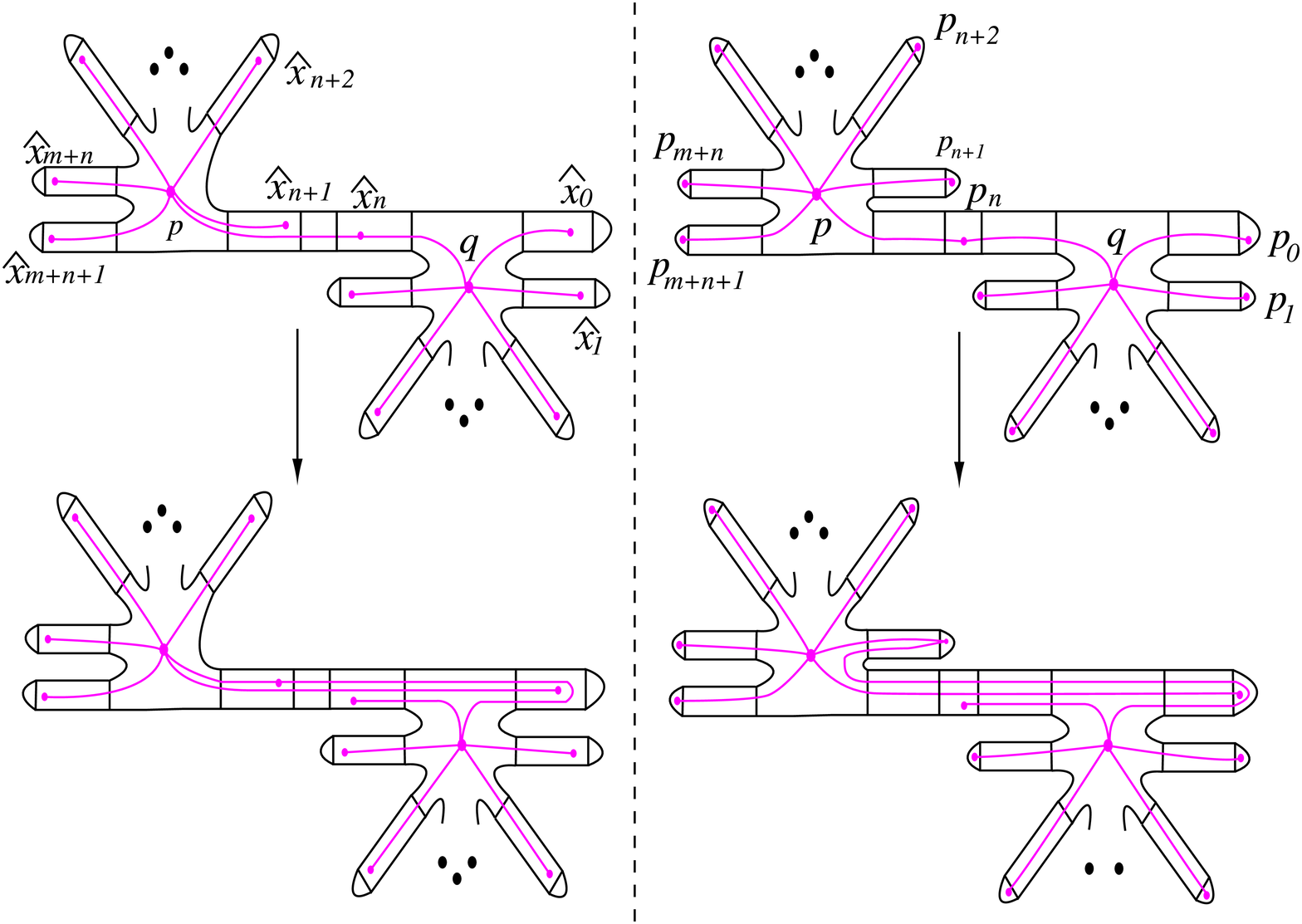}
\caption{Image of $\widehat{G}$ under $g_{m,n}$ (left), 
image of $G_{m+1,n}^{\mathrm{red}}$ under $g_{m+1,n}$ (right). }
\label{fig_proof_image}
\end{center}
\end{figure}

By Lemma~\ref{lem_symbol2} we have $[\beta_{m,n}] \succeq bt(\mathcal{O}_{g_{m,n}}(\widehat{x_0}), g_{m,n})$.  
For the proof of (1), we will show that 
$[\beta_{m+1,n}] = bt(\mathcal{O}_{g_{m,n}}(\widehat{x_0}), g_{m,n})$. 
Now we consider $\mathfrak{g}_{m+1,n}^{\mathrm{red}}$ and the thick graph map $g_{m+1,n}$ 
associated to $\mathfrak{g}_{m+1,n}^{\mathrm{red}}$ 
with the representative orbit $P_{m+1,n}$.  
Since $[\beta_{m+1,n}]=bt(P_{m+1,n}, g_{m+1,n})$, 
it suffices to show that there exists an orientation preserving homeomorphism 
$$j: (D, \mathcal{O}_{g_{m,n}} (\widehat{x_0})) \rightarrow (D, P_{m+1, n})$$
such that 
$g_{m+1,n}: D  \rightarrow D $ is isotopic to 
$j \circ g_{m,n} \circ  j^{-1}: D \rightarrow D $ relative to $P_{m+1,n}$. 
To do so, we take the tree $\widehat{G}$ embedded on 
${\Bbb F}(G_{m,n}^{\mathrm{red}}) \subset D$ (as in Figure~\ref{fig_proof_isotopy}(left)) with the following conditions: 
\begin{enumerate}
\item
$\mathcal{V}(\widehat{G})$ consists of $p,q \in \mathcal{V}(G_{m,n}^{\mathrm{red}})$ 
and all points of $\mathcal{O}_{g_{m,n}}(\widehat{x_0})$.  
\item
The valences of $p$, $q$, $\widehat{x_n} \in \mathcal{V}(\widehat{G})$ are $m+2$, $n+1$, $2$ respectively, and the other vertices have the valence $1$. 
\item
The $n+1$ edges emanate from $q$ to each $\widehat{x_0}, \cdots,  \widehat{x_n}$, and the $m+2$ edges emanate from $p$ to each 
$\widehat{x_n}, \cdots, \widehat{x_{m+n+1}}$. 
\item 
Each edge of $\widehat{G}$ transverses to each arc of ${\Bbb F}(G_{m,n}^{\mathrm{red}})$. 
\item 
$e(p, \widehat{x_{n}})$ is below $e(p, \widehat{x_{n+1}})$
with respect to the vertical coordinate of ${\Bbb F}(G_{m,n}^{\mathrm{red}})$. 
\end{enumerate}
We write $P_{m+1,n}= \{p_0, \cdots, p_{m+n+1}\}$. 
Without loss of generality we set 
$$ \mathcal{V}(G_{m+1,n}^{\mathrm{red}})= \{p, q, p_0, \cdots, p_{m+n+1}\}.$$ 
Now we take a homeomorphism $j: D \rightarrow D$ 
with $j(\widehat{G}) = G_{m+1,n}^{\mathrm{red}}$ so that 
$e(q, \widehat{x_j})$ $(0 \le j \le n)$ and 
$e(p, \widehat{x_k})$ $(n \le k \le m+n+1)$ of $\widehat{G}$ 
map to $e(q, p_j)$ and $e(p, p_k)$ of $G_{m,n}^{\mathrm{red}}$ respectively 
(Figure~\ref{fig_proof_isotopy}). 
Consider the image of $\widehat{G}$ under $g_{m,n}$ and that of $G_{m+1,n}^{\mathrm{red}}$ 
under $g_{m+1,n}$ (Figure~\ref{fig_proof_image}). 
Then  $g_{m+1,n}(G_{m+1,n}^{\mathrm{red}})$ is isotopic to 
$j \circ g_{m,n} \circ j^{-1} (G_{m+1,n}^{\mathrm{red}})$ 
relative to $\mathcal{V}(G_{m+1,n}^{\mathrm{red}})$ as union of edges. 
This means that 
$g_{m+1,n}: D \rightarrow D $ is isotopic to $j  \circ g_{m,n} \circ j^{-1}: D \rightarrow D $ 
relative to $\mathcal{V}(G_{m+1,n}^{\mathrm{red}})$. 
In particular $g_{m+1,n}: D \rightarrow D$ is isotopic to 
$j \circ g_{m,n}\circ  j^{-1}$ relative to $P_{m+1,n}$ 
since $P_{m+1,n} \subset \mathcal{V}(G_{m+1,n}^{\mathrm{red}})$. 
This completes the proof of (1). 

\begin{figure}[htbp]
\begin{center}
\includegraphics[width=4in]{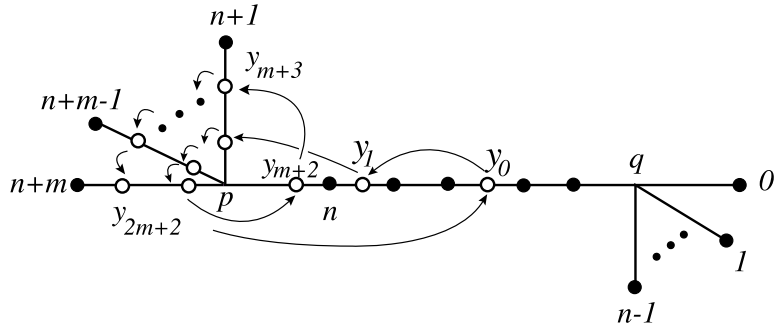}
\caption{Symbol $\circ$ indicates the periodic orbit of $y_0$, and 
$\bullet$ is a vertex of $(G_{m,n}^{\mathrm{red}})_s$. 
Note that $y_0 \in e(q,n)^3$, $y_1 \in e(q,n)^5$. }
\label{fig_smn-orbit}
\end{center}
\end{figure}

Proof of (2).  
First we show that $[\beta_{m,n}] \succeq [\sigma_{m,m+2}]$. 
We continue to consider the transition graph $\Xi_{\mathfrak{g}_{m,n}}$. 
Let  $E_0 = e(q,n)^3$, $E_1 = e(q,n)^5$, $E_2 = e(p,n+1)$, 
$ E_3= e(p,n+2) , \cdots, E_m=e(p,n+m-1), 
E_{m+1}^i=e(p,n+m)^i$, $i = 1,2$, and $E_{m+2} = e(p,n)$. 
Take the closed path $\mathcal{D}_{m+2}$ of $\Xi_{\mathfrak{g}_{m,n}}$ of length $2m+3$ 
such that 
$$E_0 \rightarrow E_1 \rightarrow \cdots \rightarrow E_m \rightarrow E_{m+1}^1 \rightarrow 
E_{m+2} \rightarrow E_2 \rightarrow E_3\rightarrow \cdots \rightarrow E_m \rightarrow 
E_{m+1}^2 \rightarrow E_0.$$
Let $y_0 \in  E_0=e(q,n)^3$ be a periodic point for $\mathfrak{g}_{m,n}^{\mathrm{red}}$ 
given in Lemma~\ref{lem_symbol}. 
Since $e(q,n)^3$ is a proper subedge of $e(q,n)$, $y_0$ is regular. 
Clearly the period of $y_0$ equals the length of $\mathcal{D}_{m+2}$, that is $2m+3$. 
Figure~\ref{fig_smn-orbit} indicates the position of the periodic orbit of $y_0$. 

Let $\mathcal{O}_{g_{m,n}}(\widehat{y_0})$ be the periodic orbit for $g_{m,n}$ 
associated to $\mathfrak{g}_{m,n}^{\mathrm{red}}$ given in Lemma~\ref{lem_symbol}. 
Recall that  $\mathfrak{h}_{m,m+2}$ is the induced graph map for $\sigma_{m,m+2}'$. 
Then we see that the braid type of $\mathcal{O}_{g_{m,n}}(\widehat{y_0})$ 
for $g_{m,n}$ equals the braid type of the representative orbit for the thick graph map 
associated to $\mathfrak{h}_{m,m+2}^{\mathrm{red}}$,  
see Figures~\ref{fig_s-graph_map}(right) regarding $n=m+2$ and Figure~\ref{fig_smn-orbit}. 
Hence  $bt(\mathcal{O}_{g_{m,n}}(\widehat{y_0}), g_{m,n})= [\sigma_{m,m+2}']= [\sigma_{m,m+2}]$.  
Since $y_0$ is regular, we obtain $[\beta_{m,n}] \succeq  [\sigma_{m,m+2}]$. 

We turn to the proof of  $[\beta_{m,n}] \succeq [\sigma_{m, m+2+\ell}]$ for any $\ell \ge 1$. 
Consider the following closed path $\mathcal{D}_{m+2+\ell}$: 
$$ \underbrace{E_0 \rightarrow E_0 \rightarrow \cdots 
\rightarrow}_{\mbox{length\ }\ell} \underbrace{E_0 \rightarrow E_1 \rightarrow \cdots \rightarrow E_{m+2} \rightarrow 
E_2 \rightarrow E_3 \rightarrow \cdots \rightarrow E_{m+1}^2 \rightarrow E_0}_
{\mbox{closed \ path\ } \mathcal{D}_{m+2}}.$$
This is the concatenation of the $\ell$--iterations of $E_0 \rightarrow E_0$ and the closed 
path $\mathcal{D}_{m+2}$. 
By using  the same argument as above, one shows that 
the braid type of the periodic orbit for $g_{m,n}$ associated to $\mathcal{D}_{m+2+\ell}$  
is $[\sigma_{m,m+2+\ell}]$, and $[\beta_{m,n}] \succeq [\sigma_{m, m+2+\ell}]$.

Proof of (3). 
Theorem 15(b) in \cite{dCH} directly shows the claim of (3). 
$\Box$ 
\medskip

Proof of Theorem~\ref{thm_horseshoe}. 
Let $Q$ be any periodic orbit in $\Omega$ for $\mathtt{H}$, and let $w_Q= (w_0 w_1 \cdots w_{s-1})$, $w_i \in \{0,1\}$ 
the code for $Q$. 
We will find a periodic orbit for the thick graph map $g_{1,k}$ associated to 
$\mathfrak{g}_{1,k}^{\mathrm{red}}$ whose braid type equals $bt(Q, \mathtt{H})$. 
We denote the edge path 
$\overline{e(p,k+1)^2} \cdot \overline{e(p,k+1)^1} \cdot  e(p,k)$ by $E_1$ and 
$\overline{e(q,k)^5} \cdot \overline{ e(q,k)^4} \cdot \overline{e(q,k)^3}$ by $E_0$. 
Then 
$\mathfrak{g}_{1,k}^{\mathrm{red}}(E_0)$ and $\mathfrak{g}_{1,k}^{\mathrm{red}}(E_1)$ 
 pass through $E_i$ $(i =0,1)$ one times, 
see Figure~\ref{fig_b-graph_map} regarding $m=1$. 
These imply that 
for the code $w_Q$, 
there exists a periodic point $z_0 \in E_{w_0}$ for $\mathfrak{g}_{1,k}^{\mathrm{red}}$ 
and a periodic point $\widehat{z_0} \in \pi^{-1}(E_{w_0})$ for $g_{1,k}$ 
such that 
\begin{eqnarray*}
(\mathfrak{g}^{\mathrm{red}}_{1,k})^s(z_0)= z_0,  &\ & 
z_i= (\mathfrak{g}_{1,k}^{\mathrm{red}})^i(z_0) \in E_{w_i} \ \mbox{and}
\\
g_{1,k}^s(\widehat{z_0})= \widehat{z_0}, &\ &  
\widehat{z_i}= g_{1,k}^i(\widehat{z_0}) \in \pi^{-1}(E_{w_i})
\end{eqnarray*}
for  each $i \in \{0, \cdots, s-1\}$. 
It is easy to check that $z_0$ is regular from the itinerary of $z_0$. 
Note that the restriction map 
$g_{1,k}|_{\pi^{-1}(E_0) \cup \pi^{-1}(E_1)}$ contracts the vertical direction of the fibered surface uniformly and expands the horizontal direction uniformly. 
Set 
$$\Omega'= \displaystyle \bigcap_{j \in {\bf Z}} g_{1,k}^j (\pi^{-1}(E_0) \cup \pi^{-1}(E_1)).$$
Then 
$g_{1,k}|_{\Omega'}: \Omega' \rightarrow \Omega'$ is conjugate to the shift map 
$\sigma: \Sigma_2 \rightarrow \Sigma_2$, and hence 
$g_{1,k}|_{\Omega'}$ is conjugate to $\mathtt{H}|_{\Omega}$. 
Thus the braid type of a periodic orbit, say $P$ in $\Omega'$ for $g_{1,k}$ 
equals that of the periodic orbit in $\Omega$ for $\mathtt{H}$ with the same itinerary as $P$. 
In particular, we have 
$bt(\mathcal{O}_{g_{1,k}}(\widehat{z_0}), g_{1,k})= bt(Q, \mathtt{H})$. 
The regularity for $z_0$ guarantees that 
$[\beta_{1,k}] \succeq bt(Q,\mathtt{H})$.  
This completes the proof. 
$\Box$
\medskip

Proof of Corollary~\ref{cor_period3}. 
By \cite[Theorem 0.2]{Han}, 
any pseudo--Anosov braid type $[mc] \in BT_3$ forces $[\sigma_1 \sigma_2^{-1}] (= [\beta_{1,1}]$). 
Since the forcing relation $\succeq$ is a partial order, 
by Theorem~\ref{thm_order1}(1) we obtain the claim of (1). 
By Theorem~\ref{thm_horseshoe} we obtain the claim of (2). 
$\Box$

\noindent
\hspace{10cm}
Eiko Kin
\\
\hspace{10cm}
Department of Mathematical and 
\\
\hspace{10cm}
Computing Sciences
\\
\hspace{10cm}
Tokyo Institute of Technology 
\\
\hspace{10cm}
Tokyo,  Japan
\\
\hspace{10cm}
e-mail: kin@is.titech.ac.jp


\end{document}